\documentclass{amsart}%\usepackage{times}
\usepackage{amssymb}

\usepackage{graphicx,epstopdf}
\begin{document}

\title
{Bernoulli numbers and solitons}%
\author{M-P. Grosset}
\address{Department of
Mathematical Sciences,
          Loughborough University,
Loughborough,
          Leicestershire, LE11 3TU, UK}
          \email{M.Grosset@lboro.ac.uk}
          
          %    Information for
%second author
         \author{A.P. Veselov}
\address{Department of
Mathematical Sciences,
          Loughborough University,
Loughborough,
          Leicestershire, LE11 3TU, UK
          and
Landau Institute for Theoretical Physics, Moscow, Russia}

\email{A.P.Veselov@lboro.ac.uk}

\begin{abstract}
We present a new formula for the Bernoulli numbers as the following integral
$$B_{2m} =\frac{(-1)^{m-1}}{2^{2m+1}} \int_{-\infty}^{+\infty} (\frac{d^{m-1}}{dx^{m-1}} \mbox{sech}^2 x )^2dx.
$$
This formula is motivated by the results of Fairlie and Veselov, who discovered the relation of Bernoulli polynomials with soliton theory.
\end{abstract}

\maketitle

\centerline{\it Dedicated to Hermann Flaschka on his
$60^{\text{th}}$ birthday}

\section{Introduction}

In the paper \cite{FV} D. Fairlie and one of the authors discovered an interesting relation of the Bernoulli polynomials with the theory of 
the Korteweg-de Vries (KdV) equation
$$u_t - 6uu_x + u_{xxx} = 0.$$
It is known since 1967 due to Gardner, Green, Kruskal and Miura \cite{MGK} that this equation has infinitely many conservation laws of the form
$$
I_m[u] = \int P_m (u, u_x, u_{xx}, ..., u_{m}) dx, $$
where $P_m$ are some polynomials of the function $u$ and its $x$-derivatives up to order $m$.
They are uniquely defined by some homogeneity property modulo adding a total derivative
and multiplication by a constant. 
This constant can be fixed by demanding that $ P_m (u, u_x, u_{xx}, ..., u_{m})  = u_{m}^2$ plus a function of derivatives of order less than $m.$

The KdV equation is famous for its remarkable family of  solutions, known as {\it solitons}, the simplest of which is a one-soliton
solution $$ u = - 2 {\rm \mbox{sech}}^2 (x - 4t), $$
corresponding to the initial profile $u(x,0) = -2 {\rm \mbox{sech}}^2 x. $

The main result of \cite{FV}  is the following formula relating the Faulhaber polynomials $F_m$
with the integrals of the KdV equation:
\begin{equation}
I_{m-1}[-2 \lambda \mbox{sech} ^2 x]= (-1)^{m-1}\frac{2^{2m+2}}{2m+1} F_{m}(\lambda).
\label{KdV_Faulhaber}
\end{equation}
Recall that the Faulhaber polynomials are directly related to the Bernoulli polynomials through the formula
\begin{equation}
B_{2m+2}(x+1) = (2m+2)F_m\left(\frac{x^2+x}{2}\right) + B_{2m+2},
\label{rel}
\end{equation}
where $B_k(x)$ and $B_{k}$ are Bernoulli polynomials and Bernoulli numbers respectively (see \cite{K, AS}).  The Bernoulli numbers have the following generating function:
$$\frac{z}{e^z-1}=\sum_{k=0}^{\infty}\frac{  B_{k}}{k!}z^k.$$
All odd Bernoulli numbers except $B_1 = - \frac{1}{2}$ are zero and the first even Bernoulli numbers are $$B_0 = 1, \, B_2 = \frac{1}{6}, \, B_4 = - \frac{1}{30}, \, B_6 =  \frac{1}{42}, \, B_8 = -  \frac{1}{30}, \, B_{10} =  \frac{5}{66}, \, B_{12}= -\frac{691}{2730},...$$ They play an important role in analysis, number theory, algebraic topology and many other areas of mathematics.

In this note we show that the relationship with soliton theory brings the following formula for the Bernoulli numbers:
\begin{equation}
\label{main-formula}
B_{2m} =\frac{(-1)^{m-1}}{2^{2m+1}} \int_{-\infty}^{+\infty} ((\mbox{sech}^2 x )^{(m-1)})^2 dx,\quad m \geq 1, 
\end{equation}
where $(\mbox{sech}^2 x )^{(m-1)}$ denotes the $(m-1)$th derivative of $\mbox{sech}^2 x.$ %\newline
We present also a direct proof of this formula based on the elegant arguments due to Logan, which we found in the book \cite{Knuth} by Graham, Knuth and Patashnik.

\section{Proofs of the main formula.}

Let us first present the proof based on the results of \cite{FV}.
Recall that the Faulhaber polynomial $F_{m}(\lambda) $ with $m \geq 1$ has the form $$F_{m}(\lambda) = \alpha_2 ^{m} \lambda ^2 + \alpha_3 ^{m} \lambda ^3+ ... +\alpha_{m+1} ^{m} \lambda ^{m+1}$$ with some rational coefficients $\alpha_2 ^{m},  \alpha_3 ^{m}, ..., \alpha_{m+1} ^{m}$ (see \cite{K}). Since the only quadratic term in the density of the KdV integral $I_{m-1}$  is $\int u_{m-1} ^2 dx,$  we have from the relation (\ref{KdV_Faulhaber}) that
\begin{equation}
\int_{-\infty}^{+\infty}  ( (\mbox{sech} ^2 x)^{(m-1)})^2dx = (-1)^{m-1}\frac{2^{2m}}{2m+1} \alpha_2 ^{m}.
\label{KdV_Faulhaber2}
\end{equation}

Differentiating the formula (\ref{rel}) twice with respect to $x$ gives
\[(\frac{2x+1}{2})^2 F_m''(\frac{x(x+1)}{2}) + F_m'(\frac{x(x+1)}{2}) =\frac{1}{2m+2}B_{2m+2}''(1+x) = (2m+1) B_{2m}(1+x)\]
because $B_{2m+2}(x)'' = (2m+2)(2m+1) B_{2m}(x).$
Since $F_m'(0)=0$ this reduces to 
\[ \frac{1}{4} F_m''(0) = ( 2m+1) B_{2m}(1).\]
Now using the well-known symmetry $B_k(1-x) = (-1)^{k} B_k (x),$ we have $B_{2m}(1)=  B_{2m}(0)= B_{2m}$ 
and thus
\[ \alpha_2 ^m  = \frac{1}{2} F_m''(0)= 2(2m+1) B_{2m}.\]
Substituting this  into (\ref{KdV_Faulhaber2}) we come to the formula (\ref{main-formula}). 

We are now going to prove the formula (\ref{main-formula}) directly without reference to soliton theory. 
We borrow the main idea from the book \cite{Knuth}, where it is attributed to Logan.

Consider the integral $J_m=\frac{(-1)^{m-1}}{2^{2m+1}}\int_{-\infty}^{+\infty} ((\mbox{sech}^2 x )^{(m-1)})^2dx.$ Integrating $J_m$ by parts $m-1$ times gives
\[J_m=\frac{1}{2^{2m+1}}\int_{-\infty}^{+\infty} (\mbox{sech}^2 x )^{(2m-2)} \mbox{sech}^2 x dx =
\frac{1}{2^{2m+1}}\int_{-\infty}^{+\infty} \tanh x ^{(2m-1)} \tanh x ^{(1)} dx.\]

Let  $y = \tanh x$ then
\[J_m=\frac{1}{2^{2m+1}}\int_{-1}^{+1} T_{2m-1}(y) dy,\]
where the polynomial $T_k(y)$ is the $k$-th derivative of $ y = \tanh x$ rewritten in terms of $y$:
\[ T_1 = y^{(1)} = 1 - y^2\]
\[ T_2 = y^{(2)} =-2 y y^{(1)}= -2 y (1 - y^2)=-2y + 2 y^3\]
\[ T_3 = y^{(3)} = -2y^{(1)} +6y^2 y^{(1)}= (1 - y^2)(-2 + 6 y^2)= -2 +8y^2 -6 y^4\]
\[ T_4 = y^{(4)} = 16 y y^{(1)} -24 y^3 y^{(1)}=(1 - y^2)(16 y - 24 y^3)=16 y - 40 y^3 + 24 y^5 \] and so on.

These polynomials can be determined by the recurrence formula 
\begin{equation}
T_n(x) = (1- x^2) T_{n-1}(x)'
\end{equation} 
with $T_0(x)=x.$ \footnote{These polynomials should not be confused with the classical Chebyshev (Tchebycheff) polynomials, which are also denoted as $T_n(x)$. We follow the notations of the book \cite{Knuth}, where trigonometric version of these polynomials was considered.} They have integer coefficients with the highest one being equal to $(-1)^n n!$ and have 
the symmetry $T_n(-x) = (-1)^{n-1}T_n(x).$ Note that the non-zero coefficients have the alternating signs
and their total sum is zero. The integers $|T_{2m-1}(0)|$ are called {\it tangent numbers} (see e.g. \cite{Knuth}). 
%By construction they satisfy the properties:
%\begin{itemize}
%\item $S_n$ is of degree $n+1$ in $x$.
%\item  $ S_n(-x) = S_n(x)$ if $n$ is odd and $ S_n(-x) = -S_n(x)$ if $n$ is even.
%\item $S_{2n +1}(0)= S_{2n }(0)' =S_{2n -1}(0)''.$
%\item $ S_n(1)=S_n(-1)= 0 $ for $n>0$. 
%\item Let $S_{n} = \sum_{k=0}^{k=n +1} a_{k }x^{k}$ then $\sum_{k=0}^{k=n +1} a_{k }=0$ and %$a_{n+1} = (-1)^ n n!.$
%\end{itemize}

Let us consider now the generating function
$$T(x,z) = \sum_{n\geq0}T_n(x) \frac{z^n}{n!}.$$

{\bf Lemma 1.} 
{\it The generating function for the polynomials $T_n(x)$ is} 
\begin{equation}
T(x,z) = \frac{\sinh z + x \cosh z}{\cosh z + x \sinh z}.
\label{generating fct}
\end{equation}

Indeed it is easy to see that  when $x = \tanh w$ the function $T(x,z)$ becomes $\tanh (z+w)$.
Now the claim follows from the Taylor formula and the definition of the polynomials $T_n(x).$

From this one can derive an interesting relation between Bernoulli and tangent numbers, see \cite{Knuth}, formula (6.93) and discussion after that.  We will need however a slightly different result.

{\bf Lemma 2.}
{\it Bernoulli numbers can be written as
\begin{equation}
B_{m} = \frac{1}{2^{m +1}} \int_{-1}^{1} T_{m-1}(x)dx  
\label{Bernoulli new def 1}
\end{equation}
for all $m>1.$}

Indeed, this formula is obvious for odd $m$ bigger than 1 since $T_{m-1}(x)$ is  an odd function and both sides of  (\ref {Bernoulli new def 1}) are then equal to $0.$ In order to prove the formula for even $m$, let us consider first the left hand side of (\ref {generating fct}). The function $T(x,z)$ can be rewritten as
\[T(x,z) = \frac{\cosh z (x+ \frac{\sinh z}{\cosh z})}{\sinh z(x+ \frac{\cosh z}{\sinh z})}  = \coth z  -\frac{1}{\sinh^2 z (x+ \coth z)}.\]
Integrating with respect to $x$ gives
\[ \int T(x,z)dx =  x\coth z -\frac{1}{\sinh^2 z} \ln |x+\coth z|, \]
which leads to
\begin{equation}
\int_{-1}^{1} T(x,z)dx = 2\coth z -\frac{2z}{\sinh^2 z} =  2\coth z  + 2z \coth' z.
\label{in}
\end{equation}
The expansion of the function $\coth z$ can be written in terms of the Bernoulli numbers:
\begin{equation}
\coth z = \frac{1}{z} + \sum_{n=1} a_{2n-1}z^{2n-1}, \quad a_{2n-1}= \frac{2^{2n}B_{2n}}{(2n)!}.
\label{co}
\end{equation}
This easily follows from the identity
$$ \frac{z}{2} \coth \frac{z}{2} = \frac{z}{2} + \frac{z}{e^z -1}.$$
From (\ref{in}) and (\ref{co}) it follows that
\begin{equation}
\int_{-1}^{1} T(x,z)dx  = 4 \sum_{n=1} n a_{2n-1} z^{2n-1}.
\label{int f{x,t} 1}
\end{equation}
On the other hand from (\ref {generating fct}) we have 
\begin{equation}
\int_{-1}^{1} T(x,z)dx  =\sum_{n\geq1}(\int_{-1}^{1}T_{2n-1}(x) dx) \frac{z^{2n-1}}{(2n-1)!},
\label{int f{x,t} 2}
\end{equation} 
since $T_{2n}(x)$ are odd polynomials.
Comparing (\ref{int f{x,t} 1}) and (\ref{int f{x,t} 2}) and replacing $a_{2n-1}$ by $\frac{2^{2n}B_{2n}}{(2n)!}$ we obtain 
\[B_{2n} = \frac{1}{2^{2n +1}} \int_{-1}^{1} T_{2n-1}(x)dx .\]
This proves Lemma 2.

{\bf Remark.} As a corollary we have the following representation of the Bernoulli number $B_{2n}$  as a sum of the fractions of the type
$$B_{2n} = \frac{1}{2^{2n}} \left( \frac{k_0}{1} + \frac{k_1}{3} + \frac{k_2}{5} + \dots + \frac{k_{n}}{2n+1}\right),$$
where the denominators are consecutive odd numbers and the numerators are the coefficients
of the polynomials $T_{2n-1}(x).$ Comparison with the formula (6.93) from the book \cite{Knuth}
leads to the following relation with the tangent numbers:
\begin{equation}
\int_{-1}^{1} T_{n}(x)dx = \frac{2n+2}{2^{n+1}-1} T_n(0).
\label{sn}
\end{equation}

% $$B_2= \frac{1}{2^{2}} (\frac{1}{1} -\frac{1}{3}), \quad B_4 = \frac{1}{2^{4}} (-\frac{2}{1} +\frac{8}{3}-%\frac{6}{5}), \quad B_6 = \frac{1}{2^{6}} (\frac{16}{1} -\frac{136}{3} +\frac{240}{5} -\frac{120}{7}).$$

Now combining Lemma 2 with the definition of the polynomials $T_k$ we have our main 

{\bf Theorem.} {\it Bernoulli numbers $B_{2m}$ with $m \geq 1$ have the following integral representation}
$$B_{2m} =\frac{(-1)^{m-1}}{2^{2m+1}} \int_{-\infty}^{+\infty} (\frac{d^{m-1}}{dx^{m-1}} \mbox{sech}^2 x )^2dx.
$$

\section{Concluding remark.}

There exists a generalisation of the Bernoulli numbers, known as {\it Bernoulli-Hurwitz numbers} $BH_k$ and related to the coefficients of the Laurent series of the Weierstrass elliptic function
$\wp(z)$ at zero:
\begin{equation}
\wp \left( z\right) =z^{-2}+\sum_{k=1}^{\infty }\frac {BH_{2k+2}}{(2k+2) (2k)!} z^{2k}
\label{expansion_of_P}
\end{equation}
(see \cite{Katz,Onishi}). 

Our formula suggests a different generalisation of Bernoulli numbers related to an elliptic curve, namely
\begin{equation}
B^{ell}_{2m} =\frac{(-1)^{m-1}}{2^{2m+1}} \int (\wp (z) ^{(m-1)})^2dz, \quad m >1, 
\label{ellip}
\end{equation}
where the integral is taking over one of the cycles on the corresponding elliptic curve.
Obviously these numbers depend also on the choice of such a cycle, which is an element of the first homology group of the curve. They naturally appear in relation with the spectral density for the classical Lam\'e operator \cite{GV}.

\section{Acknowledgements.}

We are grateful to Victor Enolski and Nick Trefethen for very useful and stimulating discussions.

%The work was supported by the Department of Mathematical Sciences of Loughborough University.

\addcontentsline{toc}{section}{\underline{References}}

 \end{document}